\documentclass[reqno,11pt]{amsart}

\newtheorem{theorem}{Theorem}[section]

\theoremstyle{definition}
\newtheorem{definition}[theorem]{Definition}

\newtheorem{proposition}[theorem]{Proposition}

\theoremstyle{remark}
\newtheorem{remark}[theorem]{Remark}

\numberwithin{equation}{section}

\usepackage[utf8]{inputenc}
\usepackage{amsmath}
\usepackage{amsfonts}
\usepackage{amssymb}
\usepackage{amsthm}
\usepackage{bbm}
\usepackage[hidelinks]{hyperref}
\usepackage{esint}

\newcommand{\eps}{\varepsilon}
\newcommand{\B}{\mathcal{B}}
\newcommand{\F}{\mathcal{F}}

\renewcommand{\H}{\mathcal{H}}

\newcommand{\R}{\mathbb{R}}
\renewcommand{\S}{\mathbb{S}}
\newcommand*{\genbf}[1]{\ifmmode\mathbf{#1}\else\textbf{#1}\fi}

\newcommand{\osc}{\mathrm{osc}}

\author{Giacomo Del Nin}
\address{Mathematics Institute, The University of Warwick, Zeeman Building, CV47HP, Coventry, UK.}
\email{\href{mailto:Giacomo.Del-Nin@warwick.ac.uk}{Giacomo.Del-Nin@warwick.ac.uk}}%

\title[Rectifiability of the jump set]{Rectifiability of the jump set of locally integrable functions}

\begin{document}

\begin{abstract}
    In this note we show that for every measurable function on $\R^n$ the set of points where the blowup exists and is not constant is $(n-1)$-rectifiable.
    In particular, for every $u\in L^1_{loc}(\R^n)$ the jump set $J_u$ is $(n-1)$-rectifiable.
\end{abstract}

\maketitle

\section{Introduction}
In the study of fine properties of functions of bounded variation a prominent role is played by the jump set, that is the set of all points where one-sided approximate limits exist (and are different) from both sides of some hyperplane. More precisely, given $u\in L^1_{loc}(\Omega)$, with $\Omega\subset\R^n$, we call $x\in \Omega$ a \emph{jump point} of $u$ if there exist $a,b\in \R$ distinct and $\nu\in\S^{n-1}$ such that
\[
\lim_{r\to 0}\fint\limits_{B_r^+(x,\nu)}|u(y)-a|dy=0\qquad\text{and}\qquad \lim_{r\to 0}\fint\limits_{B_r^-(x,\nu)}|u(y)-b|dy=0
\]
where 
\[
\begin{cases}
B_r^{+}(x,\nu)=\{y\in B_r(x): \nu\cdot (y-x) >0\}\\
B_r^{-}(x,\nu)=\{y\in B_r(x): \nu\cdot (y-x) <0\}
\end{cases}
\]
and $B_r(x)$ is the ball of center $x$ and radius $r$ (we also set $B_1:=B_1(0)$). The set of all jump points is called \emph{jump set} and denoted by $J_u$.

Among various fine properties of functions whose distributional gradient is in some form a measure, such as $BV$ or $BD$\footnote{$(G)BV$ and $(G)BD$ are respectively the spaces of functions with (generalised) bounded variation and (generalised) bounded deformation. For precise definitions we refer respectively to \cite{AFP} and \cite{ACD,Dal13}. $\H^k$ stands for the $k$-dimensional Hausdorff measure.}, it is often mentioned that the jump set $J_u$ is $(\H^{n-1},n-1)$-rectifiable, that is, it can be covered by countably many $(n-1)$-dimensional Lipschitz graphs up to an $\H^{n-1}$-negligible set \cite{AFP,ACD,Dal13,ArrSko}. We here show that this is the case for every $u\in L^1_{loc}(\R^n)$, without any requirement on the derivative of $u$, and that there is no need for the additional negligible set.
More generally, we prove the rectifiability of the set of all points where a blowup exists and is not a constant function. We call a function $v\in L^1(B_1)$ the \emph{blowup} of $u$ at the point $x\in\Omega$ if 
\[
u^{x,r}\to v\quad\text{in $L^1(B_1)$ as $r\to 0$}
\]
where $u^{x,r}\in L^1(B_1)$ are defined by $u^{x,r}(y):=u(x+ry)$\footnote{Since $\Omega$ is open, $u^{x,r}$ is well-defined for sufficiently small $r$. We could have also used the $L^1_{loc}$ convergence but we prefer to restrict to $B_1$.}. We also define the maps $T_{x,r}(y):=\tfrac{1}{r}(y-x)$, so that $u(y)=u^{x,r}(T_{x,r}(y))$. We can equivalently define $J_u$ as the set of all points where the blowup coincides with the jump function
\[
u_{a,b,\nu}(y)=\begin{cases} a & \text{if $\nu\cdot y >0$}\\
b & \text{if $\nu\cdot y <0$}
\end{cases}
\]
for some distinct $a,b\in\R$ and $\nu\in\S^{n-1}$, that is
\[
J_u=\left\{x\in \Omega: u^{x,r}\to u_{a,b,\nu}\text{ in $L^1(B_1)$ as $r\to 0$, $a\neq b\in\R$, $\nu\in\S^{n-1}$}\right\}.
\]

\begin{definition}\label{def:sing}
Given $u\in L^1_{loc}(\Omega)$ we define its \emph{singular set} as
\[
\Sigma_u:=\left\{x\in\Omega: u^{x,r}\to v \text{ in $L^1(B_1)$ as $r\to 0$, for some $v$ not constant}\right\}.
\]
\end{definition}
This set is made of all points where, provided that the blowup exists, $u$ is not approximately continuous.

\begin{theorem}\label{thm:genjump}
For every $u\in L^1_{loc}(\Omega)$,  $\Sigma_u$ (and thus also $J_u$) can be covered by countably many $(n-1)$-dimensional Lipschitz graphs.
\end{theorem}

The same result holds for vector-valued functions, working on components. In this case we just need the existence of a non-constant blowup on one component, while in the other components the limit may not even exist. The necessity of the non-constancy of the blowup is clear when we consider that a continuous function has constant blowups everywhere.

The conclusion of Theorem \ref{thm:genjump} holds even for larger sets, although we restricted the attention to $\Sigma_u$ for simplicity. We could ask for instance that the blowup exists up to adding a constant: consider the set of all points $x\in\Omega$ for which for every  $r>0$ there exist real numbers $c_r$ such that
\[
u^{x,r}-c_r\to v \text{ in $L^1(B_1)$ as $r\to 0$, for some $v$ not constant}.
\]
Then this set clearly contains $\Sigma_u$ and it is $(n-1)$-rectifiable. This fact and Theorem \ref{thm:genjump} follow from Proposition \ref{prop:osc} below which is stated in terms of the local oscillations of $u$.

\begin{definition}[Oscillation]
Given $u\in L^1_{loc}(\Omega)$ and given any Borel set $A\Subset \Omega$ with $|A|:=\H^n(A)>0$, we define the \emph{oscillation} of $u$ on $A$ as
\[
\osc(u,A):=\inf_{c \in \R} \fint_A |u(y)-c|dy.
\]
\end{definition}
Whenever it is well-defined, $\osc(u,A)=0$ if and only if $u$ is constant on $A$. For this reason the singular set $\smash{\Sigma_u}$ coincides with the set of all $x\in\Omega$ such that the blowup $v$ of $u$ at $x$ satisfies $\osc(v,B_1)>0$.

We now define a space of functions that on some ball $B\subset B_1$ oscillate less than on the whole ball $B_1$. Given $\delta>0$, $0<\tau<1$ and a ball $B\subset B_1$ we thus define
\[
\F_{\delta,\tau, B}:=\{w\in L^1(B_1): \osc(w,B_1)\geq\delta,\, \osc(w,B)\leq \tau\delta\}.
\]

\begin{proposition}\label{prop:osc}
Fix $\delta>0$, $0<\tau<1$, a ball $B\subset B_1$ and $u\in L^1_{loc}(\Omega)$. Then the set
\[
E_{\delta,\tau,B}:=\left\{x\in\Omega:\exists\, r_0>0\text{ such that } u^{x,r}\in \F_{\delta,\tau, B}\text{ for every $0<r\leq r_0$}\right\}
\]
can be covered by countably many $(n-1)$-dimensional Lispchitz graphs.
\end{proposition}

\begin{remark}
The rectifiability of the jump set is usually stated in some restricted space (e.g. in $BV$, $BD$, $GBV$, $GBD$) because it follows from the rectifiability of some other set containing it, which is harder to prove and really requires to use some information on the derivatives of $u$. For instance:
\begin{enumerate}
\item  In the space of $BV$ functions the \emph{approximate discontinuity set} $S_u$ (defined as the set of all non-Lebesgue points of $u$) is $(n-1)$-rectifiable \cite[Theorem~3.78]{AFP}. Since $J_u\subseteq S_u$, the rectifiability of $J_u$ follows from that of $S_u$, which is proved with a careful use of the coarea formula. We note also that $J_u\subseteq \Sigma_u\subseteq S_u$. 
\item In the space of $BD$ functions the rectifiability of $J_u$ follows from that of the larger set $\Theta_u$, made of all points where the upper Hausdorff $(n-1)$-density of the symmetrized gradient is positive.
In this case the rectifiability is proved thanks to a slicing formula and Besicovitch-Federer's projection theorem \cite[Proposition~3.5]{ACD}.
\item In the space $GBD$ the rectifiability of $J_u$ follows again from that of (a suitable variation of) $\Theta_u$, which is proved with arguments similar to the $BD$ case \cite[Proposition~6.1]{Dal13}.
\item In the space of functions with bounded $\B$-variation, for elliptic operators $\B$, the rectifiability of $J_u$ has been proven using rectifiability criteria for measures $\mu$ satisfying some density assumptions, and for some $\mathbb{C}$-elliptic operators $\B$ the rectifiability of a direct analogue of $\Theta_u$ holds as well \cite{ArrSko}. 
\end{enumerate}
We were quite surprised to discover that rectifiability holds for any $u\in L^1_{loc}$ and that, up to our knowledge, no mention of this was present in the literature. Since in the above cases much more than just the rectifiability of the jump set is proved, it is possible that for this reason the rectifiability of $J_u$ alone has been overlooked so far and has not been established on its own. Although we use a fairly standard decomposition argument to obtain the existence of tangent cones, we thought it worthwhile to write it down in this note to clearly separate the rectifiability of $J_u$ (and even of $\Sigma_u$) from any other further assumption on $u$.
\end{remark}

\begin{remark}
The reason why the rectifiability of $\Sigma_u$ holds without further requirements is due to the high rigidity entailed by the definition of $\Sigma_u$, in particular the existence of the blowup in the limit as $r\to 0$ and not just along some sequence $r_j\to 0$. 
In this respect we would like to mention the following result proved by Mattila \cite{Mat05}: if a measure $\mu$ on $\R^n$ has a unique tangent measure\footnote{Tangent measures are an analogue of the blowup in the space of measures, but limits are taken along any sequence $r_j\to 0$. By compactness, the existence of the limit as $r\to 0$ is equivalent to the uniqueness of all the possible limits along sequences $r_j\to 0$.} $\mu$-a.e., then $\mu$-a.e. the tangent measures must be flat, i.e. of the form $\H^k\llcorner V$ for some $k\in\{0,\ldots,n\}$ and some $k$-plane $V$. In particular, if we exclude the points where the tangent measures are multiples of $\H^n$ (the equivalent of having constant blowup), then $\mu$ has tangents of dimension at most $n-1$. 
\end{remark}

\subsection*{Acknowledgments} I would like to thank Adolfo Arroyo-Rabasa and Anna Skorobogatova for many discussions about fine properties of functions, from which this note originated, and Andrea Merlo for bringing to my attention reference \cite{Mat05}. I would also like to thank the anonymous referee for suggesting the 
simpler approach for the extension to measurable functions in Subsection \ref{sec:meas}. This project has received funding from the European Research Council (ERC) under the European Union’s Horizon 2020 research and innovation programme under grant agreement No 757254 (SINGULARITY).

\section{Proofs}

\subsection{Simplified version} 
We start with the following extremely simplified one-dimensional version of Theorem \ref{thm:genjump}, to convey the idea that will be employed later.
Consider a function $u:\R\to \R$ and, just for this paragraph, define its jump set $J_u$ as the set of all points $x$ for which there exists $r>0$ and $a,b\in \R$ distinct such that
\[
u\equiv a\quad \text{in $(x-r,x)$}\qquad\text{and}\qquad u\equiv b\quad \text{in $(x,x+r)$}.
\]
Then $J_u$ is a discrete set (and thus at most countable, i.e. $0$-rectifiable). Indeed if for $x\in J_u$ the above holds, there can be no other point of $J_u$ in $(x-r,x)\cup (x,x+r)$, because $u$ is locally constant there and thus can not jump around any point.

In the proof of the general case the local constancy is replaced by the approximate continuity of the blowup function $v$, while the presence of the jump is replaced by the non-constancy of $v$. In higher dimensions a similar argument shows that for (a suitable countable decomposition of) $J_u$ there are some directions along which there are no other points of the set, thereby proving the existence of tangent cones and thus rectifiability.

\subsection{Proof of Proposition \ref{prop:osc}}
Given a ball $B$ we denote by $r(B)$ its radius. 
\par 1. We define $E_{\delta,\tau,B,r_0}$ to be the set of all $x\in \Omega$ such that for every $0<r\leq r_0$
\begin{align}
    &\osc(u,x+rB_1)=\osc(u^{x,r},B_1 )\geq \delta \label{eq:oscB1}\\
    &\osc(u,x+rB)=\osc(u^{x,r},B) \leq \tau\delta.\label{eq:oscB}
\end{align}
Then 
    \[
    E_{\delta,\tau,B}=\bigcup_{r_0>0}E_{\delta,\tau,B,r_0}
    \]
    and the union can be taken among countably many values of $r_0$. It is thus sufficient to prove the conclusion for each of the sets $E_{\delta,B,\tau,r_0}$, which we will consider fixed from now on.
\par 2. Write $B=B_\rho(z_0)$ for some $\rho<r_0$ and $z_0\in B_1$. From the general property
\begin{equation}\label{eq:oscincrease}
\osc(w,B')\leq \frac{|B|}{|B'|}\osc(w,B)\qquad\text{if $B'\subset B$}
\end{equation}
and from \eqref{eq:oscB} we obtain that for every $0<r\leq r_0$
\begin{equation}\label{eq:osc}
\osc(u^{x,r}, B') <\delta
\end{equation}
whenever $B'$ is a ball contained in $B$ with $r(B')< \tau^{1/n}r(B)=\tau^{1/n}\rho$. In particular this is true for every ball $B'$ with radius $\rho':=\tfrac12\tau^{1/n}\rho$ and center lying in $B_\eps(z_0)$, where $\eps=\rho-\rho'$.

\par 3. We define $C$ as the convex hull of $B_\eps(z_0)\cup \{0\}$. We claim that
\begin{equation}\label{eq:cone}
E_{\delta,\tau,B,r_0}\cap (x+r_0C)=\{x\}.
\end{equation}
Indeed, suppose by contradiction this is not the case, and that there exists $x'\in E_{\delta,\tau,B,r_0}$ lying in $(x+r_0C)\setminus \{x\}$. By definition this means that there exists $z\in B_\eps(z_0)$ such that $x'-x=rz$ for some $0<r\leq r_0$. Let us consider the ball $B_{\rho'r}(x')$.
On one hand, using the map $T_{x',\rho'r}$ the ball $B_{\rho'r}(x')$ corresponds to $B_1$, and by \eqref{eq:oscB1} we have $\osc(u,B_{\rho'r}(x'))\geq\delta$. On the other hand, using the map $T_{x,r}$ the ball $B_{\rho'r}(x')$ corresponds to the ball $B_{\rho'}(z)$, and thus by \eqref{eq:osc} it satisfies $\osc(u,B_{\rho'r}(x'))<\delta$. This gives a contradiction and proves \eqref{eq:cone}.

\par 4. Property \eqref{eq:cone} says that $E_{\delta,\tau,B,r_0}$ has a one-sided tangent cone at $x$ at scale $r_0$, but the same property holds for the two-sided cone obtained symmetrizing $C$ with respect to the origin, just using the symmetry property 
\[
x'\in (x-r_0 C)\iff x\in (x'+r_0 C). 
\]
This shows that $E_{\delta,\tau,B,r_0}$ has non-trivial tangent cones at every point and thus by a standard argument (see e.g. \cite[Lemma~15.13]{Mat95}) it is contained in countably many $(n-1)$-dimensional Lipschitz graphs.\qed

\subsection{Proof of Theorem \ref{thm:genjump}}
We consider the countable family of \emph{rational balls} (i.e. with rational center and radius) that are contained in $B_1$. We pick any point $x\in \Sigma_u$ where the blowup of $u$ is a function $v$ with $\osc(v,B_1)>\delta$ for some $\delta>0$. By Lebesgue density theorem $v$ is approximately continuous at almost every point in $B_1$. In particular there exists a ball $B\subset B_1$ such that $\osc(v,B)<\tfrac12\delta$, and using \eqref{eq:oscincrease} we can find a rational ball that satisfies the same inequality. Since $v$ is the $L^1$-limit of $u^{x,r}$ we obtain the same inequalities for $u^{x,r}$, for $r$ small enough. We thus conclude that every $x\in \Sigma_u$ belongs to some set $E_{\delta, \frac12,B}$ as defined in Proposition \ref{prop:osc}, for some $\delta>0$ and some rational ball $B$. Using all rational balls and a countable sequence of $\delta$'s going to zero, we thus cover $\Sigma_u$ by countably many sets of the form $E_{\delta,\frac12, B}$, and the conclusion follows by Proposition \ref{prop:osc}.\qed

\section{Final remarks}
We chose to write the proof for the case of locally integrable functions in Euclidean spaces, but similar results hold true with very similar proofs in slightly different settings.
\subsection{Measurable functions}\label{sec:meas}
Theorem \ref{thm:genjump} extends to measurable functions possibly assuming values $\pm\infty$, even on a positive measure subset.
In this case limits have to be considered not in the $L^1$ topology but in measure. One possible approach is to fix a bounded, continuous and strictly increasing function $\Phi:\R\to\R$ (say, $\Phi(x)=\arctan x$), and then declare that $w_r\to v$ in $B_1$ as $r\to 0$ if and only if for every $\eps>0$
\[
\lim_{r\to 0}|\{x\in B_1: |\Phi(w_r(x))-\Phi(v(x))|>\eps \}|=0.
\]
This convergence agrees with the standard convergence in measure for functions that are finite a.e., and is implied by the $L^1$ convergence, but extends also to the case of functions assuming values $\pm\infty$ (it is sufficient to define $\Phi$ at $\pm\infty$ by its limits). If we define $\widetilde \Sigma_u$ as in Definition \ref{def:sing}, but replacing convergence in $L^1$ with the above convergence, then the conclusion of Theorem \ref{thm:genjump} is true for $\widetilde\Sigma_u$, for any measurable function $u$.
To prove this case it is sufficient to apply Theorem \ref{thm:genjump} to the function $\Phi\circ u$, which is now bounded and thus belongs to $L^1_{loc}(\R^n)$. The result follows considering that $\widetilde\Sigma_u=\Sigma_{\Phi\circ u}$. 

\subsection{Metric groups with dilations}
The same proof goes through, with virtually no modifications, in the setting of locally compact Lie groups with dilations. This class includes for instance all Carnot groups. 

Following \cite{Mat05} we consider a locally compact metric group $(G,\cdot)$ with a left-invariant metric $d$, i.e. such that for every $x,y,p\in G$
\[
d(p\cdot x,p\cdot y)=d(x,y).
\]
We also suppose that there is a family of dilations $\delta_\lambda: G\to G$, for $\lambda\geq 0$, that satisfy 
\[
d(\delta_\lambda x,\delta_\lambda y)=\lambda d(x,y)\quad\text{for $x,y\in G$ and $\lambda>0$}
\] 
and
\[
\delta_0=0,\qquad \delta_1=Id,\qquad\delta_{\lambda_1}\circ \delta_{\lambda_2}=\delta_{\lambda_1 \lambda_2}.
\]
We consider the Haar measure $\mu$ associated to the Lie group $G$. Given a function $u\in L^1(G,\mu)$, for every $x\in G$ and $r>0$ we define the functions $u^{x,r}\in L^1(B_1,\mu)$ by
\[
u^{x,r}(y)=u(x\cdot \delta_r y).
\] 
We then define the singular set of $u$ as
\[
\Sigma_u:=\left\{x\in G: u^{x,r}\to v \text{ in $L^1(B_1, \mu)$ as $r\to 0$, for some $v$ not constant}\right\}.
\]
Following \cite{DLMV} we give the following definition, which in the Euclidean space is equivalent to being contained in an $(n-1)$-dimensional Lipschitz graph.
\begin{definition}[Cone property]
We say that a subset $E\subseteq G$ has the cone property if there exists a ball $B\subset B_1$ (not containing the origin) such that, defining the cone $C=\{\delta_\lambda z:\lambda\geq 0, z\in B\}$, we have
\[
E\cap ( x\cdot C)=\{x\}\quad\text{for every $x\in E$}.
\] 
\end{definition}
The following version of Theorem \ref{thm:genjump} holds true in this setting.
\begin{theorem}\label{thm:genjumpgroups}
For every $u\in L^1_{loc}(G,\mu)$,  $\Sigma_u$ can be decomposed in countably many pieces having the cone property.
\end{theorem}
This theorem is proved thanks to a direct modification of Proposition \ref{prop:osc}, using the oscillations defined by
\[
\osc(u,A):=\inf_{c \in \R} \frac{1}{\mu(A)}\int_A |u(y)-c|d\mu(y).
\]
We remark that the cone property is weaker than other rectifiability notions used in Carnot groups, and we refer to \cite{DLMV} for a discussion on the subject.

\end{document}